\documentclass[10pt]{article}
\usepackage{amsmath,amsthm,amsfonts,amscd}
\DeclareMathOperator{\Hom}{Hom}
\DeclareMathOperator{\Com}{Com}
\DeclareMathOperator{\CCup}{Cup}
\DeclareMathOperator{\II}{I}
\DeclareMathOperator{\1}{id}
\DeclareMathOperator{\ad}{ad}
\DeclareMathOperator{\dev}{dev}
\DeclareMathOperator{\Ker}{Ker}
\DeclareMathOperator{\IM}{Im}
\newcommand{\NN}{\mathbb{N}}
\newcommand{\EEnd}{\mathcal End}
\newcommand{\EE}{\mathcal E}
\newcommand{\bul}{\bullet}
\newcommand{\de}{\delta}
\renewcommand{\u}{\smile}
\renewcommand{\=}{:=}
\renewcommand{\t}{\otimes}
\renewcommand{\o}{\circ}

\newtheorem{thm}{Theorem}[section]
\theoremstyle{definition}
 \newtheorem{defn}[thm]{Definition}
\theoremstyle{definition}
 \newtheorem{exam}[thm]{Example}
\begin{document}
\title{\Large\bf OPERADS AND COHOMOLOGY}
\date{}
\author{\normalsize L. Kluge and E. Paal\\ \\
\small Department of Mathematics, Tallinn Technical University\\
\small Ehitajate tee 5, 19086 Tallinn, Estonia}
\maketitle
\thispagestyle{empty}
\begin{abstract}
It is clarified how cohomologies and Gerstenhaber algebras can
be associated with linear pre-operads (comp algebras). Their
relation to mechanics and operadic physics is concisely discussed.

\end{abstract}

\section{Introduction and outline of the paper}

Operads, in essence, were invented by Gerstenhaber
\cite{Ger,Ger68} and Stasheff \cite{Sta63}. The notion of an
operad was formalized by May \cite{May72} as a tool for iterated
loop spaces. In 1994/95 \cite{GeVo94,VoGe}, Gerstenhaber and
Voronov published main principles of the operad calculus. Quite a
remarkable research activity on operad theory and its applications
can be observed in the last decade (e.~g. \cite{Rene,Smi01}). It
may be said that operads are also becoming an interesting and
important mathematical tool for QFT and deformation quantization.

In this paper, the essential parts of the operad algebra are
presented, which are relevant to understand how the cohomology
and Gerstenhaber algebra can be associated with a pre-operad. We
start from simple axioms. Basic algebraic constructions
associated with a linear pre-operad are introduced. Their
properties and the first derivation deviations of the
pre-coboundary operator are explicitly given. Under certain
condition (formal associativity constraint), the Gerstenhaber
algebra structure appears in the associated cohomology. At last,
it is also concisely discussed how operads and Gerstenhaber
algebras are related to mechanics and operadic physics.

\section{Pre-operad (composition system)}

Let $K$ be a unital associative commutative ring, and let $C^n$
($n\in\NN$) be unital $K$-modules. For \emph{homogeneous} $f\in
C^n$, we refer to $n$ as the \emph{degree} of $f$ and often write
(when it does not cause confusion) $f$ instead of $\deg f$. For
example, $(-1)^f\=(-1)^n$, $C^f\=C^n$ and $\o_f\=\o_n$. Also, it
is convenient to use the \emph{reduced} degree $|f|\=n-1$.
Throughout this paper, we assume that $\t\=\t_K$.

\begin{defn}
A linear \emph{pre-operad} (\emph{composition system}) with
coefficients in $K$ is a sequence $C\=\{C^n\}_{n\in\NN}$ of unital
$K$-modules (an $\NN$-graded $K$-module), such that the following
conditions hold.
\begin{enumerate}
\item[(1)]
For $0\leq i\leq m-1$ there exist \emph{partial compositions}
\[
  \o_i\in\Hom(C^m\t C^n,C^{m+n-1}),\qquad |\o_i|=0.
\]
\item[(2)]
For all $h\t f\t g\in C^h\t C^f\t C^g$,
the \emph{composition (associativity) relations} hold,
\[
(h\o_i f)\o_j g=
\begin{cases}
    (-1)^{|f||g|} (h\o_j g)\o_{i+|g|}f
                       &\text{if $0\leq j\leq i-1$},\\
    h\o_i(f\o_{j-i}g)  &\text{if $i\leq j\leq i+|f|$},\\
    (-1)^{|f||g|}(h\o_{j-|f|}g)\o_i f
                       &\text{if $i+f\leq j\leq|h|+|f|$}.
\end{cases}
\]
\item[(3)]
There exists a unit $\II\in C^1$ such that
\[
\II\o_0 f=f=f\o_i \II,\qquad 0\leq i\leq |f|.
\]
\end{enumerate}
\end{defn}

In the 2nd item, the \emph{first} and \emph{third} parts of the
defining relations turn out to be equivalent.

\begin{exam}[endomorphism pre-operad {\rm \cite{Ger,Ger68}}]
\label{HG} Let $A$ be a unital $K$-module and
$\EE_A^n\={\EEnd}_A^n\=\Hom(A^{\t n},A)$. Define the partial compositions
for $f\t g\in\EE_A^f\t\EE_A^g$ as
\[
f\o_i g\=(-1)^{i|g|}f\o(\1_A^{\t i}\t g\t\1_A^{\t(|f|-i)}),
         \qquad 0\leq i\leq |f|.
\]
Then $\EE_A\=\{\EE_A^n\}_{n\in\NN}$ is a pre-operad
(with the unit $\1_A\in\EE_A^1$) called the \emph{endomorphism pre-operad}
of $A$.
\end{exam}

\section{Associated operations}

Throughout this paper fix $\mu\in C^2$.

\begin{defn}
The \emph{cup-multiplication} $\u\:C^f\t C^g\to C^{f+g}$ is defined
by
\[
f\u g\=(-1)^f(\mu\o_0 f)\o_f g\in C^{f+g},
\qquad|\smile|=1.
\]
The pair $\CCup C\=\{C,\u\}$ is called a $\u$-algebra (cup-algebra) of $C$.
\end{defn}

\begin{exam}
For the endomorphism pre-operad (Example \ref{HG}) $\EE_A$ one has
\[
f\u g=(-1)^{fg}\mu\o(f\t g),
      \qquad \mu\t f\t g\in\EE_A^2\t\EE_A^f\t\EE_A^g.
\]
\end{exam}

\begin{defn}
The \emph{total composition} $\bul\:C^f\t C^g\to C^{f+|g|}$ is defined by
\[
f\bul g\=\sum_{i=0}^{|f|}f\o_i g\in C^{f+|g|},
     \qquad |\bul|=0.
\]
The pair $\Com C\=\{C,\bul\}$ is called the \emph{composition algebra} of $C$.
\end{defn}

\begin{defn}[tribraces and tetrabraces]
Define the Gerstenhaber \emph{tribraces} $\{\cdot,\cdot,\cdot\}$
as a double sum
\[
\{h,f,g\}\=\sum_{i=0}^{|h|-1}\sum_{i+f}^{|f|+|h|}(h\o_i f)\o_j g\in C^{h+|f|+|g|},
    \quad|\{\cdot,\cdot,\cdot\}|=0.
\]
The \emph{tetrabraces} $\{\cdot,\cdot,\cdot,\cdot\}$ are defined by
\[
\{h,f,g,b\}\=\sum_{i=0}^{|h|-2}\sum_{j=i+f}^{|h|+|f|-1}\sum_{k=j+g}^{|h|+|f|+|g|}
((h\o_{i}f)\o_{j}g)\o_{k}b\in C^{h+|f|+|g|+|b|}.
\]
\end{defn}

It turns out that $f\u g=(-1)^f\{\mu,f,g\}$. In general, $\CCup C$ is a
\emph{non-associative} algebra. By denoting $\mu^{2}\=\mu\bul\mu$ it turns
out that the associator in $\CCup C$ reads
\[
(f\smile g)\smile h-f\smile(g\smile h)=\{\mu^{2},f,g,h\}.
\]
Thus the \emph{formal associator} $\mu^{2}$ is an
obstruction to associativity of $\CCup C$. For the endomorphism pre-operad
$\EE_A$, $\mu^{2}$ reads as an associator as well:
\[
\mu^{2}=\mu\o(\mu\t\1_A-\1_A\t\mu),\qquad\mu\in\EE_A^2.
\]

\section{Identities}

In a pre-operad $C$, the Getzler identity
\[
(h,f,g) \=(h\bul f)\bul g-h\bul(f\bul g)
         =\{h,f,g\}+(-1)^{|f||g|}\{h,g,f\}
\]
holds, which easily implies the Gerstenhaber identity
\[
(h,f,g)=(-1)^{|f||g|}(h,g,f).
\]
The \emph{commutator} $[\cdot,\cdot]$ is defined in $\Com C$ by
\[
[f,g]\=f\bul g-(-1)^{|f||g|}g\bul f=-(-1)^{|f||g|}[g,f],\qquad|[\cdot,\cdot]|=0. \tag{G1}
\]
The \emph{commutator algebra} of $\Com C$ is denoted as
$\Com^{-}\!C\=\{C,[\cdot,\cdot]\}$. By using the Gerstenhaber identity, one
can prove that $\Com^-\!C$ is a \emph{graded Lie algebra}. The Jacobi
identity reads
\[
(-1)^{|f||h|}[[f,g],h]+(-1)^{|g||f|}[[g,h],f]+(-1)^{|h||g|}[[h,f],g]=0. \tag{G2}
\]

\section{Pre-coboundary operator}

In a pre-operad $C$, define a \emph{pre-coboundary} operator $\de\=\de_\mu$ by
\begin{align*}
-\de f&\=\ad_\mu^{right}f\=[f,\mu]\=f\bul\mu-(-1)^{|f|}\mu\bul f \\
      &\,\,=f\u\II+f\bul\mu+(-1)^{|f|}\,\II\u f,\qquad \deg\de=+1=|\de|.
\end{align*}
It turns out that $\de^{2}_\mu=-\de_{\mu^{2}}$. It follows from the Jacobi
identity in $\Com^{-}\!C$ that $\de$ is a (right) derivation of
$\Com^{-}\!C$,
\[
\de[f,g]=(-1)^{|g|}[\de f,g]+[f,\de g].
\]
But $\delta$ need not be a derivation of $\CCup C$, and $\mu^{2}$ again appears
as an obstruction:
\[
\de(f\smile g)-f\smile\delta g-(-1)^{g}\delta f\smile g
=(-1)^{|g|}\{\mu^{2},f,g\}.
\]

\section{Derivation deviations}

The \emph{derivation deviation} of $\de$ over $\bul$ is defined by
\[
\dev_\bul\de(f\t g)
   \=\de(f\bul g)-f\bul\de g-(-1)^{|g|}\de f\bul g.
\]

\begin{thm}[{\rm\cite{GeVo94,KP}}]
\label{first}
In a pre-operad $C$, one has
\[
(-1)^{|g|}\dev_\bul\de(f\t g)=f\u g-(-1)^{fg}g\u f.
\]
\end{thm}
The derivation deviation of $\de$ over $\{\cdot,\cdot,\cdot\}$ is defined by
\begin{align*}
\dev_{\{\cdot,\cdot,\cdot\}}\de\,(h\t f\t g)
    \=\de\{h,f,g\}
     &-\{h,f,\de g\}\\
     &-(-1)^{|g|}\{h,\de f,g\}-(-1)^{|g|+|f|}\{\de h,f,g\}.
\end{align*}

\begin{thm}[{\rm\cite{GeVo94,KP}}]
\label{second}
In a pre-operad $C$, one has
\[
(-1)^{|g|}\dev_{\{\cdot,\cdot,\cdot\}}\de\,(h\t f\t g)=
    (h\bul f)\u g+(-1)^{|h|f}f\u(h\bul g)-h\bul(f\u g).
\]
\end{thm}
Thus the \emph{left} translations in $\Com C$ are not
derivations of $\CCup C$, the corresponding deviations are related to
$\dev_{\{\cdot,\cdot,\cdot\}}\de$. It turns out that the \emph{right}
translations in $\Com C$ are derivations of $\CCup C$,
\[
(f\u g)\bul h=f\u(g\bul h)+(-1)^{|h|g}(f\bul h)\u g.
\]
By combining this formula with the one from Theorem \ref{second} we obtain

\begin{thm}
\label{second*}
In a pre-operad $C$, one has
\[
(-1)^{|g|}\dev_{\{\cdot,\cdot,\cdot\}}\de\,(h\t f\t g)=
    [h,f]\u g+(-1)^{|h|f}f\u[h,g]-[h,f\u g].
\]
\end{thm}

\section{Associated cohomology and Gerstenhaber algebra}

Now, it can be clarified how the Gerstenhaber algebra can be associated with a
linear pre-operad. If (formal associativity) $\mu^{2}=0$ holds, then
$\de^{2}=0$, which in turn implies $\IM\de\subseteq\Ker\de$. Then one can
form an associated cohomology ($\NN$-graded module) $H(C)\=\Ker\de/\IM\de$ with
homogeneous components
\[
H^{n}(C)\=\Ker(C^{n}\stackrel{\de}{\rightarrow}C^{n+1})/
           \IM(C^{n-1}\stackrel{\de}{\rightarrow}C^{n}),
\]
where, by convention, $\IM(C^{-1}\stackrel{\de}{\rightarrow}C^{0})\=0$.
Also, in this ($\mu^{2}=0$) case, $\CCup C$ is \emph{associative},
\[
(f\u g)\u h=f\u(g\u h), \tag{G3}
\]
and $\de$ is a \emph{derivation} of $\CCup C$. Recall from above that
$\Com^{-}\!C$ is a graded Lie algebra and $\de$ is a derivation of
$\Com^{-}\!C$. Due to the derivation properties of $\de$, the
multiplications $[\cdot,\cdot]$ and $\smile$ induce corresponding (factor)
multiplications on $H(C)$, which we denote by the same symbols. Then
$\{H(C),[\cdot,\cdot]\}$ is a \emph{graded Lie algebra}. It follows from
Theorem \ref{first} that the induced $\smile$-multiplication on $H(C)$ is
\emph{graded commutative},
\[
f\smile g=(-1)^{fg}g\smile f \tag{G4}
\]
for all $f\t g\in H^{f}(C)\t H^{g}(C)$, hence $\{H(C),\smile\}$ is an
\emph{associative graded commutative} algebra.
It follows from Theorem \ref{second*} that the  \emph{graded Leibniz rule}
holds,
\[
[h,f\u g]=[h,f]\u g+(-1)^{|h|f}f\u[h,g] \tag{G5}
\]
for all $h\t f\t g\in H^{h}(C)\t H^{f}(C)\t H^{g}(C)$. At last, it is also
relevant to note that
\[
0=|[\cdot,\cdot]|\neq|\smile|=1. \tag{G6}
\]
In this way, the triple $\{H(C),\smile,[\cdot,\cdot]\}$ turns out
to be a \emph{Gerstenhaber algebra} \cite{Ger68,GGS92Am,CGS93}.
The defining identities of the Gerstenhaber algebra are (G1)-(G6).

In the case of an endomorphism pre-operad, the Gerstenhaber algebra
structure appears on the Hochschild cohomology of an associative algebra
\cite{Ger}.

\section{Discussion: x-mechanics}

Some people like commutative diagrams. Consider the following one:
\[
\begin{CD}
\text{Poisson algebras}@<\text{algebraic abstraction}<<\text{mechanics}\\
      \wr                       @.                        \wr\\
\text{Gerstenhaber algebras}@<\text{algebraic abstraction}<<\text{x-mechanics}
\end{CD}
\]
\emph{Poisson algebras} can be seen as an algebraic abstraction
of mechanics.
Here $\sim$ means \emph{similarity}:
Gerstenhaber algebras are graded analogs of the Poisson algebras.

It may be expected that there exists a kind of mechanics (x-mechanics)
associated with operads and Gerstenhaber algebras.
According to the diagram, x-mechanics is a graded analogue of mechanics and
\emph{observables} of an x-mechanical model must satisfy the
(homotopy \cite{GeVo94,VoGe}) Gerstenhaber algebra identities.

Cohomologies and Gerstenhaber algebras associated with pre-operads are natural
objects for modelling x-mechanical systems.
Physically relevant examples of the Gerstenhaber algebras and \emph{odd symplectic}
structures are provided by the Batalin-Vilkovisky algebras \cite{Ge,Sch,Ik}.
Relevance of the operad structure in handling the renormalization problems in QFT
were recently stressed in \cite{Kr,LM}.

\section*{Acknowledgement}
Research was in part supported by the Estonian Science Foundation
Grant 5634.

\end{document}